\documentclass{gtart}

\def\ifplaintex{\expandafter\ifx\csname documentclass\endcsname\relax}

\def\gtp{{\mathsurround=0pt\it $\cal G\mskip-2mu$eometry \&\ 
$\cal T\!\!$opology $\cal P\!$ublications}}  

\def\recd{{\small Received:\qua\receiveddate\ifx\reviseddate\relax
\else\qquad Revised:\qua\reviseddate\fi\par}} 

\def\lognumber#1{\def\thelognumber{#1}}
\def\volumenumber#1{\def\thevolumenumber{#1}}
\def\volumeyear#1{\def\thevolumeyear{#1}}
\def\papernumber#1{\def\thepapernumber{#1}}
\def\pagenumbers#1#2{\def\startpage{#1}\def\finishpage{#2}}
\def\published#1{\def\publishdate{#1}}

\def\received#1{\def\receiveddate{#1}}

\def\accepted#1{\def\accepteddate{#1}}

\def\asciiaddress#1{\def\theasciiaddress{#1}}

\def\asciikeywords#1{\def\theasciikeywords{#1}}

\let\\\par\let\thelognumber\relax\let\thevolumenumber\relax
\let\thepapernumber\relax\let\thevolumeyear\relax\let\startpage\relax
\let\finishpage\relax\let\publishdate\relax\let\receiveddate\relax
\let\reviseddate\relax\let\accepteddate\relax\let\theasciititle\relax
\let\theasciiauthors\relax\let\theasciiaddress\relax
\let\theasciiabstract\relax\let\theasciikeywords\relax

\let\theasciiemail\relax

\ifplaintex
\font\logobig=cmssbx10 scaled 3836
\font\logomed=cmssbx10 scaled 2557
\else
\font\logobig=cmssbx10 scaled 4200
\font\logomed=cmssbx10 scaled 2800
\fi

\long\def\makeagttitle{   
\count0=\startpage
\agt\hfill      
\hbox to 45truept{\vbox to 0pt{\vglue -13truept{\logomed A\kern -.37em{\logobig 
T}\kern -.38em G}\vss}\hss}
\break
{\small Volume \thevolumenumber\ (\thevolumeyear)
\startpage--\finishpage\nl
Published: \publishdate}

\vglue .25truein

{\parskip=0pt\leftskip 0pt plus
1fil\def\\{\par\smallskip}{\Large\bf\thetitle}\par\medskip} \vglue
0.05truein

{\parskip=0pt\leftskip 0pt plus 1fil\def\\{\par}{\sc\theauthors}
\par\medskip}%
 
\vglue 0.03truein 

{\small\leftskip 25truept\rightskip 25truept{\bf Abstract}\stdspace\theabstract

{\bf AMS Classification}\stdspace\theprimaryclass
\ifx\thesecondaryclass\relax\else; \thesecondaryclass\fi\par
{\bf Keywords}\stdspace \thekeywords\par}\vglue 7truept

}   

\ifplaintex
\font\phead=cmsl9 scaled 950
\font\pnum=cmbx10 scaled 913
\font\pfoot=cmsl9 scaled 950
\headline{\vbox to 0pt{\vskip -4.5mm\line{\small\phead\ifnum
\count0=\startpage ISSN 1472-2739 (on-line) 1472-2747 (printed)
\hfill {\pnum\folio}\else\ifodd\count0\def\\{ }%
\ifx\theshorttitle\relax\thetitle\else\theshorttitle\fi\hfill{\pnum\folio}
\else\def\\{ and }{\pnum\folio}\hfill\ifx\theshortauthors\relax\theauthors
\else\theshortauthors\fi\fi\fi}\vss}}
\footline{\vbox to 0pt{\vglue 0mm\line{\small\pfoot\ifnum\count0=\startpage
\copyright\ \gtp\hfill\else
\agt, Volume \thevolumenumber\ (\thevolumeyear)\hfill\fi}\vss}}
\else
\font=cmsl9 scaled 1050
\font\lnum=cmbx10 
\font=cmsl9 scaled 1050
\makeatletter
\def\@oddhead{{\small\lhead\ifnum\count0=\startpage ISSN 1472-2739 
(on-line) 1472-2747 (printed)\hfill {\lnum\number\count0}\else\ifodd\count0
\def\\{ }\ifx\theshorttitle\relax \thetitle \else\theshorttitle\fi\hfill
{\lnum\number\count0}\else\def\\{ and }{\lnum\number\count0}
\hfill\ifx\theshortauthors\relax 
\theauthors\else\theshortauthors\fi\fi\fi}}\def\@evenhead{\@oddhead}
\def\@oddfoot{\small\lfoot\ifnum\count0=\startpage\copyright\ \gtp\hfill\else
\agt, Volume \thevolumenumber\ (\thevolumeyear)\hfill\fi}
\def\@evenfoot{\@oddfoot}
\makeatother
\fi
\let\maketitlepage\makeagttitle
\let\makeshorttitle\maketitlepage
\let\maketitle\maketitlepage

\newwrite\gtoutfile
\long\gdef\makeheadfile{  
{\def\\{, }\def\s{ }
\immediate\openout\gtoutfile head.xxx
\immediate\write\gtoutfile{To: math@arxiv.org}
\immediate\write\gtoutfile{Subject: put OR rep NNNNN:ppppp}
\immediate\write\gtoutfile{--text follows this line--}
\immediate\write\gtoutfile{Proxy-for: \ifx\theasciiauthors\relax
\theauthors\else\theasciiauthors\fi\s<\ifx\theasciiemail\relax\theemail\else\theasciiemail\fi>}
\immediate\write\gtoutfile{\noexpand\\}
\immediate\write\gtoutfile{Authors: \ifx\theasciiauthors\relax
\theauthors\else\theasciiauthors\fi}
{\def\\{ }\immediate\write\gtoutfile{Title: \ifx\theasciititle\relax
\thetitle\else\theasciititle\fi}}
\immediate\write\gtoutfile{Subj-class: GT or SG, GR etc}
\immediate\write\gtoutfile{MSC-class: \theprimaryclass\ifx\thesecondaryclass\relax\else, \thesecondaryclass\fi}
\immediate\write\gtoutfile{Journal-ref: Algebr. Geom. Topol. \thevolumenumber\s
(\thevolumeyear) \startpage-\finishpage}
\immediate\write\gtoutfile{Comments: Published by Algebraic and
Geometric Topology at}
\immediate\write\gtoutfile{\s\s\s  http://www.maths.warwick.ac.uk/agt/AGTVol\thevolumenumber/agt-\thevolumenumber-\thepapernumber.abs.html}
\immediate\write\gtoutfile{\noexpand\\}
\immediate\write\gtoutfile{}
\ifx\theasciiabstract\relax
\immediate\write\gtoutfile{\theabstract}\else
\immediate\write\gtoutfile{\theasciiabstract}\fi
\immediate\write\gtoutfile{}
\immediate\write\gtoutfile{\noexpand\\}
\immediate\write\gtoutfile{}
\immediate\closeout\gtoutfile}}  

\def\maketitlepage{\makeagttitle\makeheadfile}
\let\makeshorttitle\maketitlepage
\let\maketitle\maketitlepage

\lognumber{18}
\volumenumber{3}
\volumeyear{2003}
\papernumber{18}
\published{19 June 2003}
\pagenumbers{557}{568}
\received{7 June 2003}
\accepted{16 June 2003}

\usepackage{amsmath,amsthm,amssymb}
\usepackage{graphicx}

\font\block=cmcsc10

\newtheorem*{thm}{Theorem}
\newtheorem*{lemma1}{Lemma 1}
\newtheorem*{lemma2}{Lemma 2}

\theoremstyle{definition}
\newtheorem*{remarks}{Remarks}
\newtheorem*{example}{Example}
\newtheorem*{remark}{Remark}
\newtheorem*{finalremarks}{Final Remarks}

\newcommand{\bz}{\mathbb Z}
\newcommand{\e}{\varepsilon}
\newcommand{\lk}{\textup{lk}}
\newcommand{\mb}{\bar\mu}
\newcommand{\ti}{{\textsc i}}
\newcommand{\tj}{{\textsc j}}

\def\place#1#2#3{\text{\kern-#1pt{\raise#2pt\hbox{\rlap{#3}}}\kern#1pt}}
\def\figure#1#2{\includegraphics[scale=#2]{#1.eps}}
\def\figurelabel#1{\centerline{\small #1}}

\begin{document}

\title{A geometric interpretation of\\Milnor's triple linking numbers}
\author{Blake Mellor\\Paul Melvin}
\address{Loyola Marymount University, One LMU Drive\\Los 
Angeles, CA 90045, USA}
\secondaddress{Bryn Mawr College, 101 N merion Ave\\Bryn Mawr,
Pa 19010-2899, USA}
\asciiaddress{Loyola Marymount University, One LMU Drive\\Los 
Angeles, CA 90045, USA\\and\\Bryn Mawr College, 101 N merion Ave\\Bryn Mawr,
Pa 19010-2899, USA}
\email{bmellor@lmu.edu, pmelvin@brynmawr.edu}

\begin{abstract}
Milnor's triple linking numbers of a link in the 3-sphere are
interpreted geometrically in terms of the pattern of intersections of the
Seifert surfaces of the components of the link.  This generalizes the well
known formula as an algebraic count of triple points when the pairwise linking
numbers vanish.
\end{abstract}

\primaryclass{57M25}
\secondaryclass{57M27}
\keywords{$\mb$-invariants, Seifert surfaces, link homotopy}
\asciikeywords{mu-bar-invariants, Seifert surfaces, link homotopy}

\maketitle


\section{Introduction}

Milnor's higher order linking numbers (or {\sl mu-bar invariants})
$\mb_\ti$ of an
oriented link
$L = L_1\cup\cdots\cup L_\ell$ in the $3$-sphere are a measure of how deep the
longitudes
$l_i$ of the components $L_i$ lie in the lower central series of the link group
$\pi = \pi_1(S^3-L)$.  They are defined as follows \cite{milnor}:

Any diagram of $L$ gives rise to a Wirtinger presentation of $\pi$ whose
generators are meridians of the link, one for each arc in the diagram.  If one
works modulo any term $\pi^n$ in the lower central series of $\pi$ (where
$\pi^1=\pi$ and $\pi^n=[\pi,\pi^{n-1}]$) then only one meridian $m_i$ is needed
from each component $L_i$.  In particular the longitude $l_i$ can be written
mod $\pi^n$ as a word $l_i^n$ in these preferred meridians, i.e.\ as an
element of the free group $F = F(m_1,\dots,m_\ell)$.  Now each element in $F$
can be viewed as a unit in the ring $A$ of power series in the noncommuting
variables $h_1,\dots,h_\ell$ by substituting $1+h_i$ for $m_i$ and
$1-h_i+h_i^2-+\cdots$ for $m_i^{-1}$; this is its {\it Magnus expansion}
\cite{magnus} (which embeds $F$ in $A$).  For any sequence $\ti=i_1 \cdots \,
i_r$ of integers between $1$ and $\ell$, let $\mu_\ti$ denote the  coefficient
of $h_{i_1} \cdots h_{i_{r-1}}$ in the Magnus expansion of $l_{i_r}^n$ for
large $n$.  It is easily shown that $\mu_\ti$ is independent of $n$ for $n\ge
r$, but unfortunately it {\sl does} depend on the original choice of
meridians.  If however one reduces modulo the greatest common divisor
$\delta_\ti$ of the lower order coefficients $\mu_\tj$ (for all proper
subsequences $\tj$ of $\ti$) then one obtains Milnor's linking numbers
$$
\mb_\ti =  \mu_\ti \ (\text{mod}\ \delta_\ti) \in
\bz/\delta_\ti\bz.
$$
They are invariants of $L$ (with a given ordering of its components) up to
link concordance, and in fact up to link homotopy if the indices
$i_1,\dots,i_r$ are distinct.  Note that in general $\mb_\ti$ depends on
the ordering of the indices in $\ti$, although it is symmetric in these
indices for $r=2$ (indeed $\mu_{ij} = \mb_{ij}$ is the usual pairwise linking
number) and antisymmetric for $r=3$.

The computation of $\mb_\ti$ from this point of view is tedious and
unenlightening, especially for large $r$.  As a warm-up consider the pairwise
linking number $\mb_{ij}$, which is the coefficient of $h_i$ in the Magnus
expansion of $l_j^2$.  To compute $\mb_{ij}$, write the  word for the longitude
$l_j$ as a product of generators in the Wirtinger presentation, each 
of which is
a conjugate $wm_pw^{-1}$ of its associated  (preferred) meridian 
$m_p$.  Working
mod the commutator subgroup $\pi^2$, this conjugate can be replaced by just
$m_p$, giving a word $u_j$ in $m_1,\dots,m_\ell$.  Then
$$
\mb_{ij}=\e_i(u_j)
$$
where $\e_i(u_j)$ is the sum of the exponents of $m_i$ in $u_j$ (cf.\
\cite[p\,8]{fenn}).

Now consider the triple linking number $\mb_{ijk}$ for distinct
$i,j,k$, which is the coefficient of $h_ih_j$ in the Magnus expansion of
$l_k^3$, reduced modulo the greatest common divisor $\delta = \delta_{ijk}$
of the pairwise linking numbers of $L_i$, $L_j$ and $L_k$.  As before, write
the word for the longitude $l_k$ as a product of generators in the Wirtinger
presentation, and then express each generator as a conjugate $wm_pw^{-1}$ of
its  associated meridian.  Working mod $\pi^3$, each generator in the
conjugating word $w$ can be replaced by its associated meridian, giving the
word $l_k^3$ in $m_1,\dots,m_\ell$.  Since we are only interested in the
coefficient  of $h_ih_j$ in the Magnus expansion of this word, any 
occurrence of
meridians  other than $m_i$ and $m_j$ can be ignored.  Thus we are left with a
word $u_k$ in $m_i$ and $m_j$.  Now
$$
\mb_{ijk} = \e_{ij}(u_k) \ (\text{mod}\ \delta)
$$
where $\e_{ij}(u_k)$ denotes the sum of the {\it signed occurences} of
$m_im_j$ in the word $u_k$: each $\cdots m_i^r\cdots m_j^s \cdots$ in $u_k$,
where $r,s=\pm1$, contributes $rs$ to $\e_{ij}(u_k)$ (see
\cite[pp\,38,142]{fenn}).  (In fact $\mb_{ijk}$ can be computed as
$\e_{ij}(u)$ for any word $u$ representing $l_k$ mod $\pi^3$.)

\pagebreak

The purpose of this note is to reinterpret this process for computing 
the triple
linking numbers in terms of Seifert surfaces of the link components and their
intersections.  This generalizes work of Cochran~\cite{cochran} for links with
trivial pairwise linking numbers (also see Turaev \cite{turaev}).

\begin{remarks}  It is worth noting a few properties of $\e_i$ and
$\e_{ij}$.  Let $u$ and $v$ be words in the meridians $m_1,\dots,m_\ell$.  Then

(1)\quad  $\e_i(uv) = \e_i(u) + \e_i(v)$

(2)\quad  $\e_{ij}(uv) = \e_{ij}(u) + \e_{ij}(v) + \e_i(u)\e_j(v)$

(3)\quad  $\e_{ij}(u) + \e_{ji}(u) = \e_i(u)\e_j(u)$

(cf.\ \cite[Lemma 4.2.7]{fenn}).  The third equation identifies two
ways of counting (with sign) all pairs of an occurrence of $m_i$ and
an occurrence of $m_j$ (in either order).
\end{remarks}

The authors gratefully acknowledge the hospitality of the Mathematics
Institute at the University of Warwick during the Geometry/Topology Workshop in
July 2000, where this work was initiated.


\section{The formula for $\mb_{ijk}$}

Fix three components $L_i,L_j,L_k$ of $L$, with $\delta$ the greatest common
divisor of their pairwise linking numbers, and choose any associated oriented
Seifert surfaces $F_i,F_j,F_k$ in general position.  Set $F=F_i\cup F_j\cup
F_k$.  It will be shown below that $\mb_{ijk}$ is the difference of two
geometric invariants, $t_{ijk}(F)\in\bz$ and
$m_{ijk}(F)\in\bz/\delta\bz$, defined
as follows.

The intersection $F_i\cap F_j\cap F_k$ consists of isolated triple points,
each of which is given a sign according to the orientations and ordering of the
surfaces: the sign is positive if and only if the ordered basis of
normal vectors
to $F_i,F_j,F_k$ at the triple point agrees with the standard orientation on
$S^3$.  Let
$$
t_{ijk}(F) = \# F_i\cap F_j\cap F_k,
$$
the total number of triple points, {\sl counted with sign}.  This integer is
clearly invariant under isotopies of the individual Seifert surfaces which
maintain their mutual general position, and is antisymmetric in the indices.


Similarly any pairwise intersection of the Seifert surfaces is a union of
circles, clasps and ribbons (see Figure 1).


\medskip{\small
\centerline{\kern30pt\figure{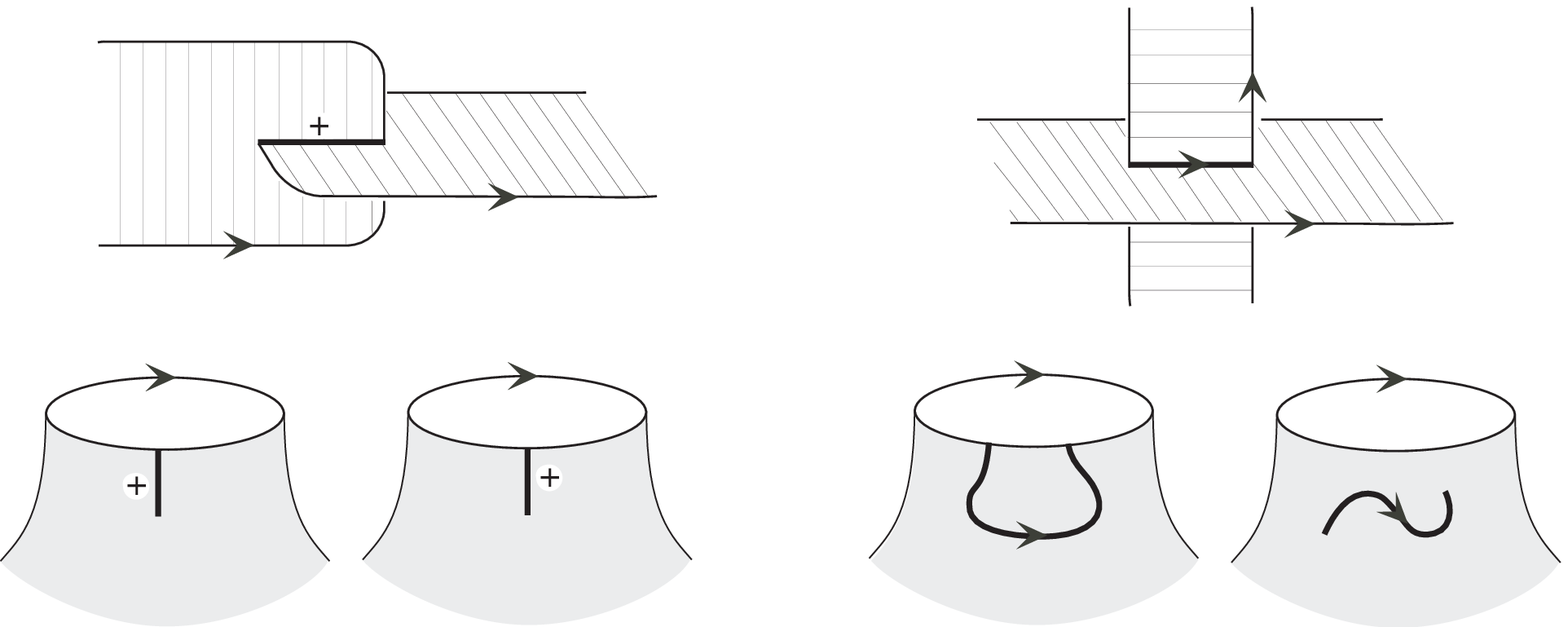}{.675}
\place{324}{77.5}{$F_i$} \place{270}{90}{$F_j$} \place{128.7}{135}{$F_i$}
\place{40.5}{108}{$F_j$}
\place{380.7}{13.5}{$F_i$} \place{247.5}{13.5}{$F_j$} \place{180}{13.5}{$F_i$}
\place{45}{13.5}{$F_j$}}
\bigskip
\centerline{(a) clasp \kern180pt (b) ribbon}
\medskip
\figurelabel{Figure 1: Intersections of Seifert surfaces}
\medskip}

\noindent
The endpoints of the clasps and ribbons are
oriented in a natural way since they are intersection points of one oriented
surface with the oriented boundary of another.  Observe that the two
endpoints of
any clasp always have the same orientation, while the endpoints of a ribbon are
oppositely oriented, and so these orientations can be specified by a
single sign
for each clasp and an orientation on each ribbon arc (from negative
to positive)
as shown in the figure.

Now choose a basepoint on each link component (this is equivalent to
choosing the
meridians in the definition of Milnor's invariant).  Starting at the basepoint
of $L_i$, read off a word $w_i$ in $j$ and $k$ by proceeding around $L_i$: each
double point in $L_i\cap F_j$ of sign $\e$ contributes $j^\e$, and
similarly for
$k$.  (Here and below we simply use the indices $i,j,k$ to denote the
corresponding
meridians $m_i,m_j,m_k$.)  Similarly define words $w_j$ and $w_k$.  Let
$e_{pqr}(F)$ denote the coefficient of $h_ph_q$ in the Magnus expansion of
$w_r$, i.e.\ $e_{pqr}(F) = \e_{pq}(w_r)$, and set
$$
m_{ijk}(F) = e_{ijk}(F) + e_{jki}(F) + e_{kij}(F),
$$
the sum of the $e_{ijk}$'s over cyclic permutation of the indices.


\begin{example}
The Borromean rings are shown in Figure 2, bounding three disks with
two clasps and
one ribbon.  If the basepoints are chosen near the labels, we compute
$w_i = 1$, $w_j = kk^{-1} = 1$, and $w_k = iji^{-1}j^{-1}$.  Thus
$e_{ijk} = 1$ and
$e_{jki} = e_{kij} = 0$, so $m_{ijk} = 1$.  Note that the clasps can be
eliminated, at the expense of introducing a triple point, by sliding the disk
$F_j$ to the left so that its boundary encircles $F_i$ and pierces $F_k$
(cf.\ the theorem below).
\end{example}


\centerline{\small\kern 45pt\figure{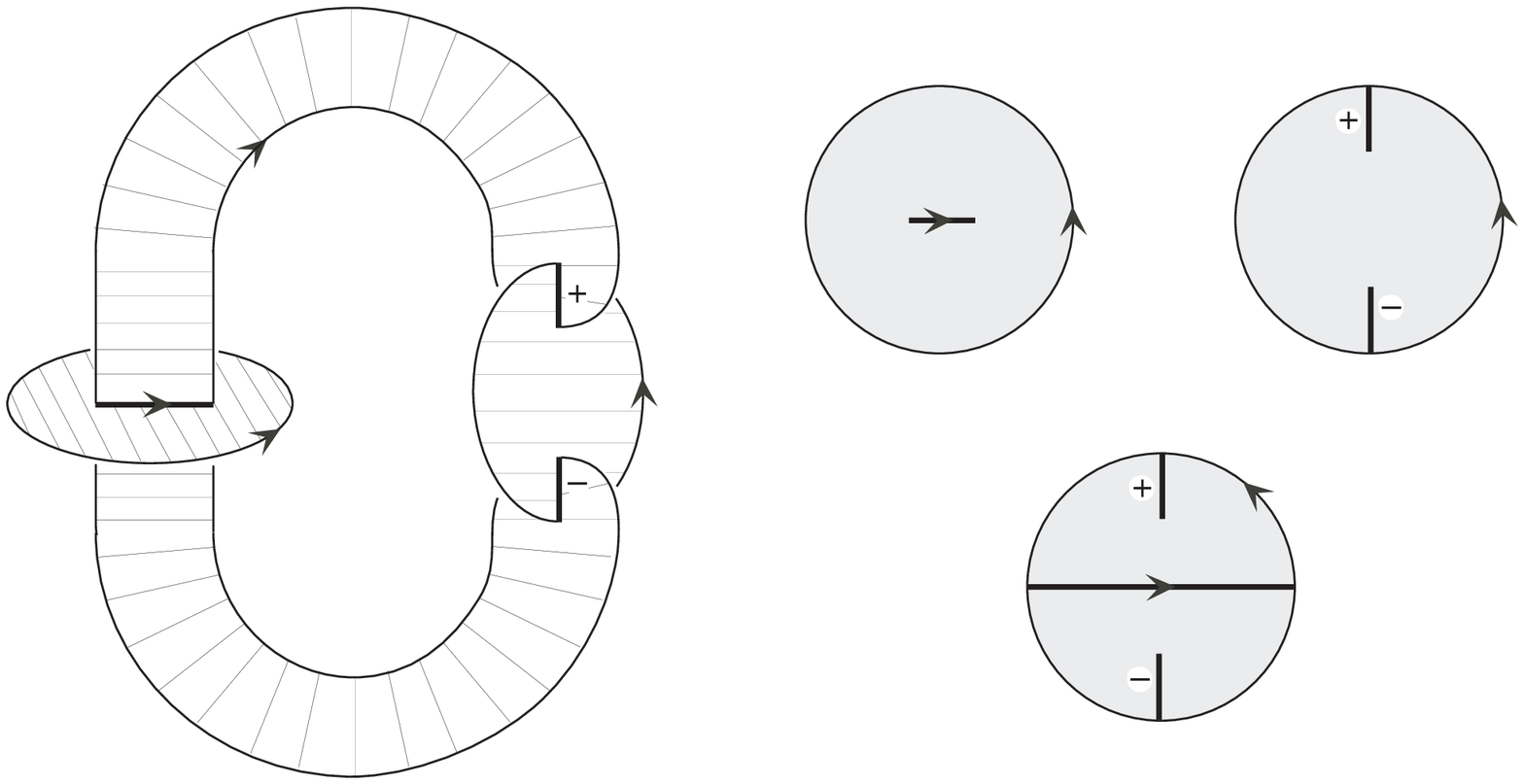}{.60}
\place{353}{85}{$F_i$} \place{195}{85}{$F_j$}
\place{175}{150}{$F_i$} \place{20}{150}{$F_j$}
\place{282}{32}{$F_k$} \place{77}{17}{$F_k$}
\place{62}{162}{$k$} \place{66}{85}{$k^{-1}$}
\place{157}{45}{$i^{-1}$} \place{80}{44}{$i$}
\place{125}{5}{$j^{-1}$} \place{124}{82}{$j$}}
\medskip
\figurelabel{Figure 2: Computation of $m_{ijk}$}

\begin{remarks} (1)\qua  Choosing a different basepoint changes
$e_{ijk}$, and hence $m_{ijk}$, by some multiple of the linking numbers, so
$m_{ijk}$ (like $\mb_{ijk}$) is well defined modulo $\delta$.  To see 
this, note
that a basepoint change induces a cyclic permutation of $w_k$, which reduces
the value of $e_{ijk}$ by $\e_{ij}(uv)-\e_{ij}(vu)$ for some factorization
$w_k=uv$.  But (using Remarks 1 and 2 at the end of the previous section)
\begin{align*}
\e_{ij}(uv) - \e_{ij}(vu) &= \e_i(u)\e_j(v) - \e_i(v)\e_j(u) \\
&= (\e_i(u)+\e_i(v))\e_j(v) - \e_i(v)(\e_j(u)+\e_j(v)) \\
&= \e_i(uv)\e_j(v)-\e_i(v)\e_j(uv) \ = \ \mu_{ik}\e_j(v) + \e_i(v)\mu_{jk}.
\end{align*}
where $\mu_{pq} = \lk(L_p,L_q)$, so this difference is trivial mod $\delta$.

\smallskip
(2)\qua   The invariants $e_{ijk}$ are antisymmetric in $i$ and $j$.  Indeed
$e_{ijk}+e_{jik} = \e_{ij}(w_k)+\e_{ji}(w_k) = \e_i(w_k)\e_j(w_k) =
\mu_{ik}\mu_{jk} \equiv 0 \pmod\delta$, by Remark 3 in \S1.  It
follows that $m_{ijk}$ is totally antisymmetric.
\end{remarks}

Our goal is to prove the following result, where as usual $\delta$ denotes the
greatest common divisor of the pairwise linking numbers of $L_i$,
$L_j$ and $L_k$.

\begin{thm}
$(1)$\qua For any choice of Seifert surfaces $F$,
$$
\mb_{ijk} \equiv m_{ijk}(F)-t_{ijk}(F) \pmod\delta.
$$

\smallskip\noindent
$(2)$\qua There is a choice of $F$ such that $t_{ijk}(F)=0$, so
$\mb_{ijk} \equiv m_{ijk}(F) \pmod\delta$.

\smallskip\noindent
$(3)$\qua There is a choice of $F$ such that $m_{ijk}(F)\equiv 0$, so
$\mb_{ijk} \equiv -t_{ijk}(F) \pmod\delta$.
\end{thm}

The familiar fact that the Milnor invariants $\mb_{ijk}$ are antisymmetric in
the indices follows immediately from the theorem and the preceding remark.  We
also recover Cochran's result that $\mb_{ijk} = -t_{ijk}(F)$ if the
boundary of $F$
is disjoint from the double point set (which is only possible if the pairwise
linking numbers vanish).

The idea of the proof is as follows.  We will define two special types of
of Seifert surfaces $F$ and show that $m_{ijk}(F) \equiv \mb_{ijk}$ and
$t_{ijk}(F)=0$ for those of type 1, while $m_{ijk}(F) \equiv 0$ for
those of type
2.  We will then describe a set of {\it finger moves} that can transform an
arbitrary $F$ into either reduced type, and that preserve $m_{ijk}-t_{ijk}
\pmod\delta$.  The theorem follows.


\section{Proof of the theorem}

A generic family $F = F_i\cup F_j\cup F_k$ of Seifert surfaces for the link
$L_i\cup L_j \cup L_k$ will be called {\it reduced} if it is of one of the
following types:

\begin{enumerate}
\item[] {\block Type 1 : (disjoint clasps)} The only double curves are disjoint
clasps.
\smallskip
\item[] {\block Type 2 : (boundary ordered)} Around the boundary of any one of
the three surfaces, all intersections with either of the other two are
adjacent, except possibly for cancelling pairs (i.e.\ adjacent intersections of
opposite sign with the second surface that may occur between intersections
with the third).
\end{enumerate}

If $F$ is a family of type 1, then $t_{ijk}(F) = 0$ (since there are no triple
points at all) and $m_{ijk}(F)$ is just the $\mb$-invariant (as is shown in the
lemma below).
If $F$ is of type 2, then by definition the word $w_k$ (see \S2) is of the form
$m_i^pm_j^q$ for an appropriate choice of basepoint on $L_k$, where $p =
\mu_{ik}$ and  $q=\mu_{jk}$ (the linking numbers of $L_i$ and $L_j$
with $L_k$) and
similarly for $w_i$ and $w_j$.  Thus $e_{ijk}(F) = \e_{ij}(w_k) = pq \equiv 0
\pmod\delta$, and so $m_{ijk}(F) \equiv 0$.  It will follow from the
theorem that
in this case $t_{ijk}(F)$ reduces (mod $\delta$) to the negative of the
$\mb$-invariant.

\begin{lemma1}
If $F$ is of type $1$, then $m_{ijk}(F) \equiv \mb_{ijk}$.
\end{lemma1}

\begin{proof}
By the classification of surfaces, each component of $F$ can be viewed as a
2-dimensional handlebody, i.e.\ a {\it disk} (0-handle) with {\it bands}
(1-handles) attached, while the clasps between any two components can be viewed
as pairs of {\it feelers}, i.e.\ bands attached at one end to their respective
disks and clasped with each other at the other end.

Of course the bands and feelers may be highly linked and twisted, so
consider the
effect of a {\it band pass} (moving one such through another or
through itself).
Clearly $m_{ijk}(F)$ is unchanged, since the pattern of double curves
is unaltered,
and $\mb_{ijk}$ is certainly unchanged under self-passes and passes
between bands
from the same surface component, since it is a link homotopy
invariant.  In fact
$\mb_{ijk}$ is also invariant under passes between bands from
distinct components.
To see this, it suffices (using the antisymmetry of $\mb$) to consider a pass
between bands on $F_j$ and $F_k$ as shown below.


\medskip
\centerline{\small\figure{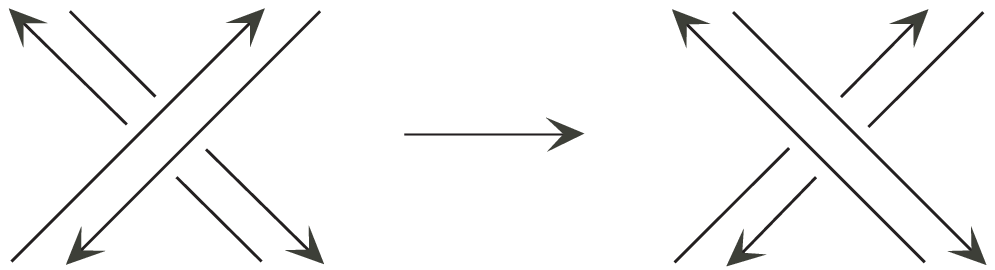}{.60} \place{178}{50}{$j$}
\place{135}{50}{$k$}
\place{70}{50}{$j$} \place{25}{50}{$k$}
\place{66}{20}{$j_1$} \place{40}{20}{$j_2$}}

\noindent
If the longitude of $L_k$ on the left is the word $u\,v$ (where $u$
starts at the
lower left point) then on the right it is $^{(j_1^{-1}j_2)}u\,v$ (where $^ab$
denotes the conjugate $aba^{-1}$) for suitable conjugates $j_1,j_2$ of the
meridian $j$ of $L_j$.  But
$$
\e_{ij}(^{(j_1^{-1}j_2)}uv) = \e_{ij}(^{(j^{-1}j)}uv) = \e_{ij}(uv),
$$
since in general $\e_{ij}(\cdots {^{(^ab)}c} \cdots) = \e_{ij}(\cdots {^bc}
\cdots)$, and so $\mb_{ijk}$ remains unchanged.

Using band passes, the three surfaces can be ``disentangled", except for the
clasps.  In particular, they can be placed so that their disks and
bands lie near
the vertices of a triangle whose edges are formed by the feelers
reaching out to
clasp each other (see Figure 3 for an example).  Furthermore, it can
be arranged
that the only link crossings aside from the clasps (under the projection to the
plane of the triangle) are between feelers from the same surface component, or
between bands from the same component.  In fact the latter can be changed at
will, since self-crossings of a component do not contribute to $m_{ijk}$ or
$\mb_{ijk}$.  Thus we may assume that each Seifert surface consists only of a
disk with some feelers,  as shown in the
figure.

Now $\mb_{ijk} = \e_{ij}(u_k)$, where $u_k$ is the word defined in \S1 which
represents the longitude $l_k$.  For the case at hand, $u_k$ is
clearly a product
of conjugates of $i$ and $j$ and their inverses.  In particular, each positive
clasp of $F_i$ with $F_k$ contributes a conjugate $j^{-p}ij^p$, where
$p$ is the
number of clasps of $F_i$ with $F_j$ (counted with sign) which
precede the given
clasp along $L_i$.  (Actually the conjugating element is a word in $i$ and
$j$ with total exponent $p$ in $j$, but the $i$'s can be ignored when computing
$\e_{ij}$.)  Likewise each negative clasp contributes $j^{-p}i^{-1}j^p$, and so
the total contribution from the clasps between $F_i$ and $F_k$ is
$e_{jki}$ (since
$\e_{ij}(j^{-p}i^{\pm1}j^p) = \pm p$).  A similar argument shows that the
contribution from the clasps between $F_j$ and $F_k$ is $-e_{ikj}$ (the minus
sign arising since $\e_{ij}(i^{-p}j^{\pm1}i^p) = \mp p$).  So the total
internal contributions from these factors is $e_{jki}-e_{ikj} =
e_{jki}+e_{kij}$ (by Remark 2 in \S2).  The contribution from the order of
the factors is $e_{ijk}$, so we conclude that $\mb_{ijk} =
e_{ijk}+e_{jki}+e_{kij} \equiv m_{ijk}$ as desired.
\end{proof}

\centerline{\small\figure{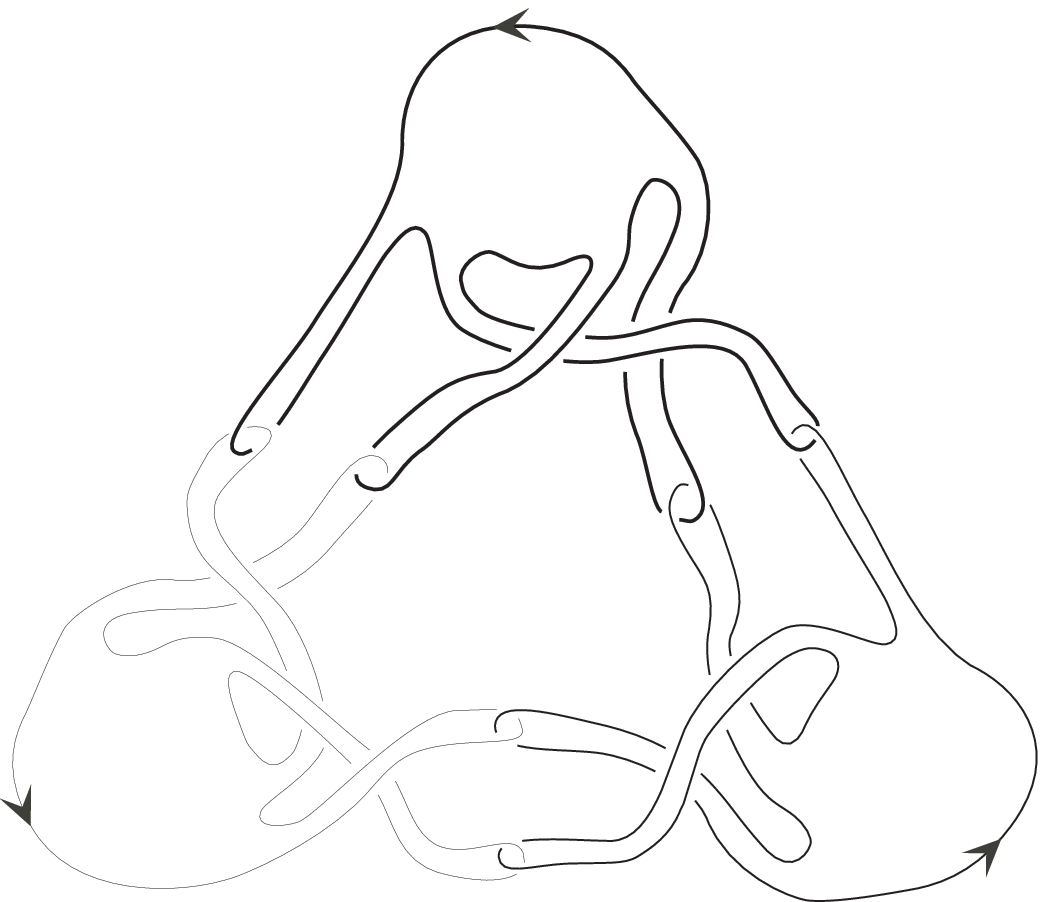}{.72}  \place{225}{4.5}{$F_i$}
\place{9}{4.5}{$F_j$} \place{90}{171}{$F_k$}}
\bigskip
\figurelabel{Figure 3: Disentangled surfaces}

\begin{remark}
There is a nice alternative proof that $\mb_{ijk} \equiv m_{ijk}(F)$ for
disentangled surfaces:  One observes that the braiding of the feelers can be
accomplished by ``Borromean tangle" replacements -- familiar from the theory
of finite type invariants of 3-manifolds -- one for each braid
generator.  One then
shows that the invariants $\mb_{ijk}$ and $m_{ijk}$ change by an equal amount
($\pm1$) under such a tangle replacement, and it is easy to check
that they agree
if there is no braiding.
\end{remark}

We now define three ``finger moves" on a generic family $F$ of
Seifert surfaces.
Each move involves pushing the interior of a small disk on one of the surfaces
along a path in one of the other surfaces.  It is understood that the
path begins
at the center of the disk, and that the center should be pushed
slightly beyond the
endpoint of the path.

\begin{enumerate}
\item[] {\block First finger move} : Push from a double point along a path to
the boundary which contains no other double points.

\smallskip
\item[] {\block Second finger move} : Push from a triple point along a path to
the boundary which lies in a double curve and contains no other triple points.

\smallskip
\item[] {\block Third finger move} : Push from a double point in the boundary
to an adjacent double point along a path in the boundary.
\end{enumerate}

\noindent
These moves are illustrated in Figures 4, 5 and 6.


\bigskip
\centerline{\small\kern 30pt\figure{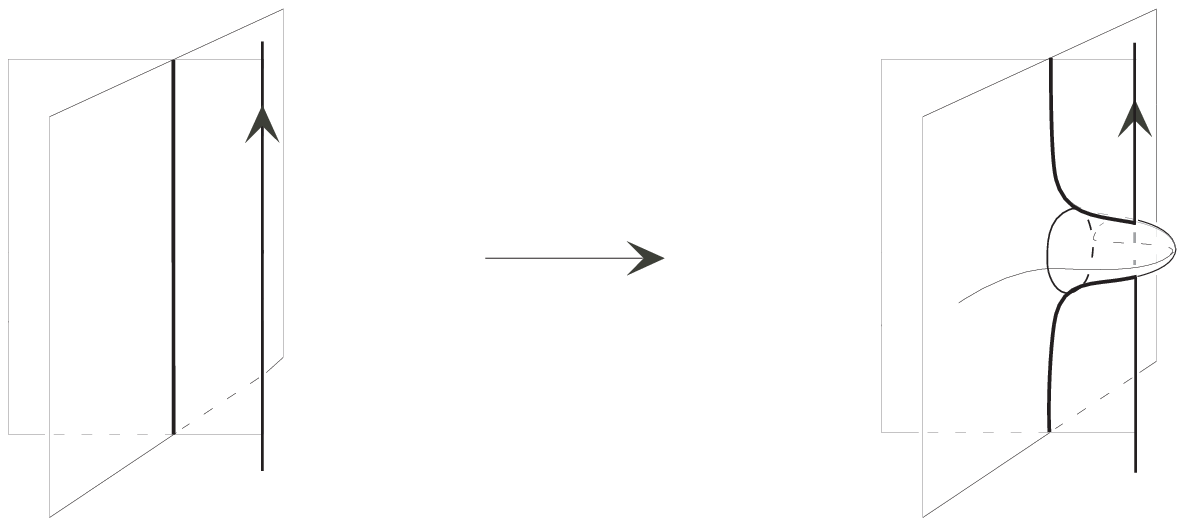}{.80}
\place{215}{10}{$L_i$} \place{295}{90}{$F_i$}
\place{290}{-5}{$F_j$}
\place{25}{10}{$L_i$} \place{107}{90}{$F_i$}
\place{100}{-5}{$F_j$}}

\bigskip
\figurelabel{Figure 4: First finger move}
\bigskip


\bigskip
\centerline{\small\kern 30pt\figure{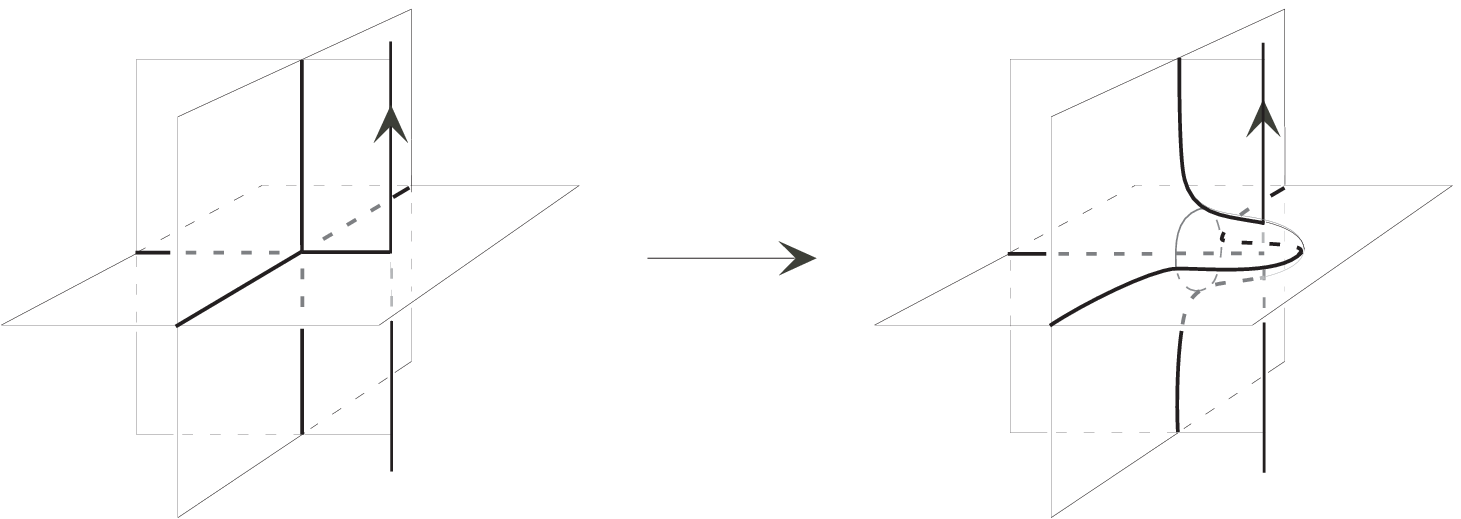}{.80}
\place{245}{10}{$L_i$} \place{325}{90}{$F_i$}
\place{320}{-5}{$F_j$} \place{230}{82}{$F_k$}
\place{57}{10}{$L_i$} \place{140}{90}{$F_i$}
\place{133}{-5}{$F_j$} \place{45}{82}{$F_k$}}

\bigskip
\figurelabel{Figure 5: Second finger move}
\bigskip


\bigskip
\centerline{\small\kern 30pt\figure{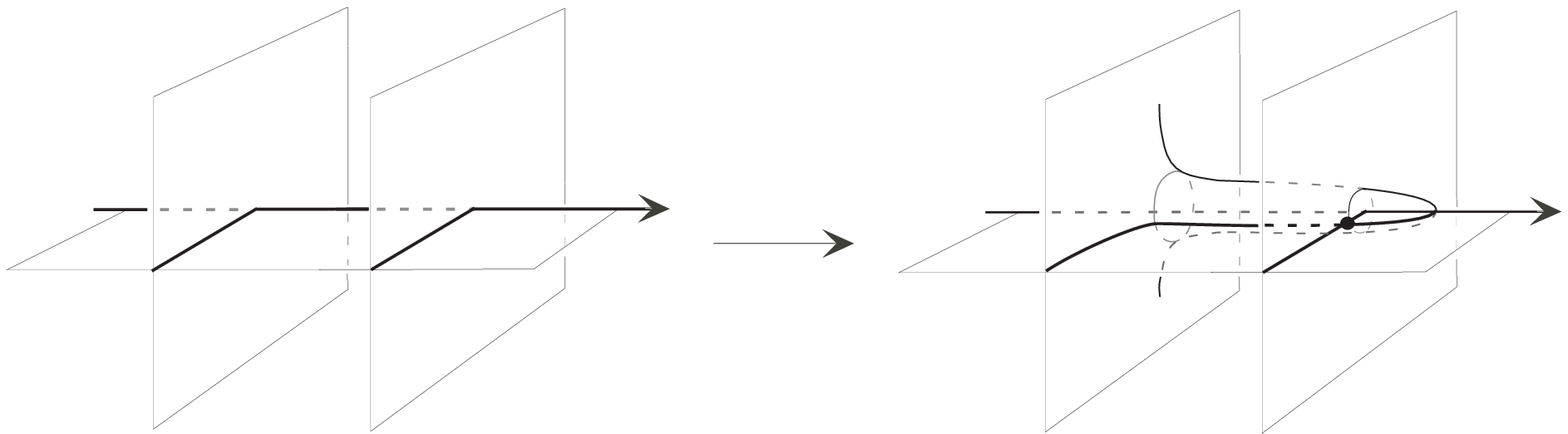}{.64}
\place{212}{56}{$L_k$} \place{352}{24}{$F_k$}
\place{320}{-8}{$F_i$} \place{272}{-8}{$F_j$}
\place{28}{56}{$L_k$} \place{168}{24}{$F_k$}
\place{136}{-8}{$F_i$} \place{88}{-8}{$F_j$}
}

\bigskip
\figurelabel{Figure 6: Third finger move}
\bigskip

\begin{lemma2}
The difference $m_{ijk}-t_{ijk}$ is invariant under all three finger moves.
\end{lemma2}

\begin{proof}
Each move is, to be exact, six possible moves, since the roles of the three
surfaces can be permuted.  We will consider only one of the six cases, the ones
depicted in the figures, since the others can be proved similarly.

The first finger move (Figure 4) clearly does not change the number of triple
points, so $t_{ijk}$ remains the same.  Also, since the words $w_j$ and $w_k$
(around the boundaries of $F_j$ and $F_k$) are unchanged, so are $e_{ijk}$ and
$e_{kij}$.  The word $w_i$ is only changed by adding a cancelling pair, so
$e_{jki}$ is unchanged.  Therefore $m_{ijk}$ is fixed as well, and so is the
difference $m_{ijk}-t_{ijk}$.

The second finger move changes the number of triple points by 1.  In the case
shown in Figure 5, a positive triple point is removed, so $t_{ijk}$ decreases
by 1 (the boundaries of the surfaces are oriented counterclockwise).  The words
$w_j$ and $w_k$ are unchanged; the word $w_i$ is modified by replacing one
appearance of the letter $k$ by its conjugate $j^{-1}kj$.  So $e_{ijk}$ and
$e_{kij}$ are unchanged, while $e_{jki}$ decreases by 1.  The result is to
decrease $m_{ijk}$ by 1.  Therefore the difference $m_{ijk}-t_{ijk}$
is unchanged.

The third finger move, as shown in Figure 6, decreases $t_{ijk}$ by 1 by
adding a negative triple point.  The words $w_i$ and $w_j$ are unchanged, but
$w_k$ is modified by replacing an adjacent pair of letters $ij$ by $ji$.  So
$e_{jki}$ and $e_{kij}$ remain fixed, while $e_{ijk}$ decreases by 1.
This means
$m_{ijk}$ is reduced by 1, and once again the difference $m_{ijk}-t_{ijk}$ is
preserved.
\end{proof}

The final step in the proof of the theorem is to show how to use finger moves
to transform a generic family of Seifert surfaces into either one of the two
reduced types.  To achieve the first type, in which the intersections
between the
components are disjoint clasps, it is necessary to eliminate all
circles, ribbons
and triple points.  We describe a procedure to accomplish this.

First, remove circles by finger moves of type 1.  Do this by following any path
from the circle to the boundary of one of the surfaces in which it
lies, breaking
each double arc crossed by this path by a type 1 finger move.  Note that this
move never creates circles, so its repeated application eventually
eliminates all
the circles.

So now all triple points involve clasps and/or ribbons.  We can inductively
remove the triple points along each clasp or ``long" ribbon (i.e.\ a double arc
with both endpoints on the boundary of one surface) by using finger
moves of type
2.  Again, this move never creates circles, so we are left with a disjoint
collection of clasps and ribbons.

Finally, remove the ribbons by finger moves of type 1 on ``short"
ribbons (double
arcs with both endpoints in the interior of the surface) along any
path from the
ribbon to the boundary which avoids other ribbons and clasps.  Such a path will
exist because there are no triple points or circles, so no short ribbon can be
separated from the boundary.  This leaves a collection of disjoint
clasps, i.e.\
a reduced family of type 1.

To obtain a reduced family of surfaces of the second type, begin with
one of the
first type, and then simply reorder the clasps as necessary via transpositions,
using the third finger move.  This completes the proof of the theorem.
\qed

\begin{finalremarks} (1)\qua The theorem suggests a natural
construction
of a three component link with any prescribed pairwise linking
numbers $p,q,r$ and
any given triple linking number $m$:  Start with the unlink, bounding disjoint
disks.  Next add $p$ feelers between one pair of disks (as in Lemma
1), and similarly for $q$ and $r$ with the other pairs.  Finally add $m$
Borromean tangles (band sums with the three components of the Borromean rings).
Now apply the theorem, using the calculation in the Example in \S2.

(2)\qua  If one allows nonorientable Seifert surfaces, then the geometric invariants
$t_{ijk}$ and $m_{ijk}$ are still defined modulo 2.  (Note that the mod
2 reduction of $\e_{ij}(w)$ depends only on the {\sl parity} of the
exponents of
$i$ and $j$ in $w$.)  Furthermore, the proof of the theorem carries
over, without
change, to show that $\mb_{ijk} \equiv m_{ijk} - t_{ijk} \pmod2$.  Of
course this
is a vacuous statement unless the pairwise linking numbers are even.

(3)\qua   A geometric interpretation for the higher order Milnor invariants remains
open.  A formula for the {\sl first non-vanishing} invariant is given by
Cochran \cite{cochran}, using surfaces associated to certain {\sl derived
links}.  It is natural to ask whether the ideas of this paper can be
extended to the general case; the difficulty seems to lie in finding the right
"higher-order surfaces."
\end{finalremarks}


\Addresses\recd

\end{document}